\documentclass[letterpaper, 10 pt, conference, english]{ieeeconf}

\usepackage{color}
\definecolor{gray}{RGB}{128,128,128}

\IEEEoverridecommandlockouts
\usepackage{epsfig}
\usepackage{subfigure}
\usepackage{amssymb,amsmath,amsfonts,layout,graphicx}
\usepackage{makeidx}
\usepackage{babel}
\usepackage{sublabel}
\usepackage{cases}
\usepackage{algorithm}
\usepackage{algorithmic}

\usepackage
[
        letterpaper,
        left=1.91cm,
        right=1.91cm,
        top=1.91cm,
        bottom=2cm,
]
{geometry}

\newtheorem{proposition}{Proposition}
\newtheorem{theorem}{Theorem}
\newtheorem{definition}{Definition}
\newtheorem{corollary}{Corollary}
\newtheorem{remark}{Remark}
\newtheorem{example}{Example}

\DeclareMathOperator*{\argmax}{arg\,max}

 \title{Transactive Energy System: Market-Based Coordination of Distributed Energy Resources}
  \author{Sen Li, Jianming Lian, Antonio Conejo, and Wei Zhang
 \thanks{This work was partly supported by the National Science Foundation under Grant CNS-1552838.}
\thanks{S. Li is with the Department Mechanical Engineering,  University of California, Berkeley, USA. (e-mail: \tt\footnotesize lisen1990@berkeley.edu)}
\thanks{J. Lian is with the Optimiza and Control Group, Pacific Northwest National Laboratory, Richland, WA 99352. (e-mail: \tt\footnotesize jianming.lian@pnnl.gov)}
\thanks{A. Conejo is with the Department of Electrical and Computer Engineering, The Ohio State University, Columbus, OH 43210. (e-mail: \tt\footnotesize conejo.1@osu.edu)}
\thanks{W. Zhang is with the Department of Mechanical Engineering,  The Southern University of Science and Technology, China. (e-mail: zhangw3@sustech.edu.cn)}
}

\begin{document}
\maketitle

\section{Introduction}
Due to pressing environmental concerns, there is a global consensus to commit to a sustainable energy future. Germany has embraced Energiewende -- a bold sustainable energy policy of no operational nuclear plants by 2022. California has set an ambitious goal that mandates $50\%$ renewable penetration by 2025, $60\%$ by 2030, and $100\%$ by 2045 \cite{sb100california}.  The vast integration of renewable energy into the power grid  imposes daunting challenges to the conventional supply-side control paradigm. First of all, renewable energy is intermittent and uncertain. Integrating renewable energy will reduce the system inertia, inject undesirable  variability, and substantially increase the need for system reserves. These reserves are typically provided by conventional generators. However, conventional generators have limited capacity and ramping rate, and they produce carbon emissions that defeat the carbon benefit of renewable integration. On the other hand, the increasing penetration of renewable energy also squeezes the response time needed to balance the power grid. That is, the time to make critical operating decisions is decreasing from minutes to seconds and in some cases even sub-seconds due to the increasing variability of supply and demand. The shorter response time for decision-making  places significant challenges on the conventional power grid that requires human interaction.

A promising pathway for large-scale renewable integration is through  distributed energy resources (DERs). DERs are those small-scale assets connected to the distribution power network \cite{DNV2014DER}, including  distributed generators, electricity storage, smart loads, among others. They are typically close to the load and can be used individually or in aggregation to reduce the power from the transmission, enhance grid security, decrease carbon emissions, and allow the deferral of new transmission lines and large power plants \cite{carr2008survey}.  Various types of control strategies have been proposed to engage both customer-owned and third-party DERs into the power system operation. The first one is direct control, where utility companies remotely control the end-use loads based on prior mutual agreements. The second one is the price control, where  pricing signals are sent to end-use loads to affect local demand. Common examples of price control include Time of Use (TOU) pricing, Critical Peak Pricing (CPP), and Real Time Pricing (RTP). Recently, transactive control is a new type of control strategy that automates and engages DERs through the market interaction. This is an emerging control paradigm originated from the Pacific Northwestern National Laboratory. Built upon the idea of price control, transactive control derives the pricing signal by coordinating individual DERs at the system level through the market. In this way, DERs are engaged in price discovery.

The GridWise Architecture Council (GWAC) first coined the term ``transactive energy" for this new concept. It is defined as ``a set of economic and control mechanism that allows the dynamic balance of supply and demand across the entire electrical infrastructure using value as a key operational parameter". Before 2015, several field demonstration projects were implemented in  the U.S. and Europe to showcase  applications of transactive control in power systems and testify its technology feasibility in real system deployment. These demonstration projects include the Olympic Peninsula Demonstration (2006--2007)~\cite{Fuller11,Hammerstrom07}, the AEP Ohio gridSMART\textsuperscript{\textregistered} Demonstration (2010--2014)~\cite{Widergren14Report,Widergren14ISGT}, the Pacific Northwest Smart Grid Demonstration (2010--2015)~\cite{Huang10,Hammerstrom15a}, and PowerMatching City (2009--2015)~\cite{Kok12}. Since then, transactive control has received substantial attentions in the research community and many research results have been published \cite{widergren2014residential}, \cite{jin2012simulation}, \cite{bejestani2014hierarchical}, \cite{junjie2016transactive}, \cite{hao2017transactive}, \cite{annaswamy2018transactive}, \cite{TEPNNLreport2018}, \cite{liu2017transactive}. Transactive energy systems are deeply intertwined with  many topics in microeconomics and control, including competitive equilibrium \cite{chen2010two}, \cite{nguyen2012walrasian}, contract theory \cite{yang2014risk}, \cite{balandat2014contract}, mechanisms design \cite{samadi2012advanced}, \cite{zhang2015truthful}, Stackelberg games\cite{maharjan2013dependable}, to name a few. It also involves numerous applications such as thermostatically controlled loads \cite{li2016market}, \cite{kalsi2012development}, \cite{hao2017transactive}, electric vehicles \cite{ma2013decentralized}, \cite{robu2011online}, \cite{liu2013planning}, \cite{li2014integrated}, \cite{liu2018two}, distributed generators \cite{rahimi2010demand}, \cite{rahbari2014incremental},  among others.


One challenge in understanding the literature is the lack of a unifying framework to systematically analyze and compare the existing results. In the literature,  many works are closely related but intricately different. If there is a slight change in one problem with respect to the information structure, the decision-making order or the rationality assumption, it could lead to another problem that is completely different in nature. Unfortunately, these intricacies are often neglected in the literature, which causes significant efforts in understanding the existing results and developing the new results. This gives rise to the following questions: how to formally compare different transactive energy systems? How to formulate a proper transactive energy system to solve a given problem? What tools are available in the literature for each class of transactive energy systems? 


This paper answers these questions by synthesizing a unifying framework that standardizes the formulation of transactive energy systems to facilitate the analysis of most transactive energy systems studied in the literature. This framework consists of four key elements: agent preference, control decision, information structure, and solution concept. Agent preference specifies the preferred outcome by each agent and models the conflict of interest in decision making. Control decision specifies the set of decisions for each agent to maximize its payoff. Information structure describes the information available to each agent before it makes a decision. Solution concept encodes the rationality assumptions of each agent in the system. These elements are important in identifying and distinguishing different transactive energy systems. In this paper, we summarize a number of important classes of transactive energy systems including the competitive market, Stackelberg game, reverse Stackelberg game, and mechanism design. Each class of transactive energy systems has a specific formulation under the proposed unifying framework. The connections and differences between different classes are discussed. Available tools and results for each class are also surveyed.


 \section{Unifying Framework for Transactive Energy Systems}
To substantiate the necessary development of a unifying mathematical framework for transactive energy systems, consider the following power system problem as a simple motivating example that well illustrates the underlying subtleness and complexity associated with the design of transactive energy systems. Suppose that one geographic area receives the electricity from the main power grid through a radial transmission line that is capacity constrained. Due to the significant growth of the population in this area, the capacity of the connection to the main grid is expected to become inadequate to supply all the cooling demand during hot summer conditions. In this case, transmission line congestion may occur when all the residents in the area  turn on their air conditioners for cooling, which increases the risk of transformer failure due to significant overloading. To address this problem, transmission upgrade is a straightforward but expensive solution. If a large amount of residential air conditioners were equipped with smart thermostat that permits advanced control, the deployment of transactive energy system can be considered as a non-wires alternative to infrastructure upgrading for congestion management. The basic idea is to treat the constrained transmission capacity as a scarce resource, and create a market among all customers within the area to realize resource allocation. The demand from individual customers will be coordinated by appropriate market signals so that the total demand from the area can be maintained below the capacity limit. The context of this problem is very straightforward. However, when it comes to the actual design of transactive energy system, many practical questions arise:
\begin{itemize}
\item What is the objective of the system coordinator?
\item What are the objectives of individual household owners who has controllable air conditioners?
\item Would market participants consider their market power on prices?
\item Would market participants share their private information with the coordinator?
\item What form of signals shall be communicated between system operator and market participants?
\end{itemize}
The answers to these questions greatly vary according to the underlying assumptions over the system. As long as the answer to any one of these questions is different, the resulting mathematical formulations could be different as well. Therefore, this necessitates a unifying framework that systematically pinpoints the conceptual difference between different formulations.

In this section, we propose such a unifying framework for market-based coordination of DERs. It consists of four elements that are essential to specify a transactive energy system. These elements include agent preference and constraint, control decision, information structure, and solution concept, which are elaborated in more details in the following subsections.

\subsection{Agent Preference}
Consider a group of $N$ DERs connected to the distribution power networks. A coordinator employs economic principles and tools to coordinate these DERs, forming a transactive energy system. We refer to this system as a multi-agent system consisting of $N+1$ agents that are indexed by $i$, where $i=0$ denotes the coordinator agent, and $i\in\{1,\ldots,N\}$ the resource agents. In this system, individual DERs have their preferences over the energy generation or consumption, which is referred to as energy allocation in the following. The coordinator have its own preference over the energy allocations for individual DERs as well. We encode these distinct, and possibly private, preferences into payoff functions. These payoff functions describe the mapping from the energy allocations and energy prices for individual DERs to real numbers that quantify the preferences of the coordinator agent and resource agents, respectively.

Let $a_i\in\mathcal{A}_i$ be the energy allocation for resource agent $i$, where $\mathcal{A}_i$ is the set of local constraints on the energy allocation $a_i$. Denote $a=(a_1,\ldots,a_N)$ and $\mathcal{A}=\mathcal{A}_1\times\cdots\times\mathcal{A}_N$. Let $\lambda_i \in \mathbb{R}$ be the energy price for resource agent $i$, where $\Lambda_i$ is the set of admissible prices for $\lambda_i$. Denote $\lambda=(\lambda_1,\ldots,\lambda_N)$ and $\Lambda=\mathbb{R}^N$. Without loss of generality, it is assumed that the payoff functions of individual DERs could depend on the energy allocations and energy prices of other DERs. Then, the payoff function of resource agent $i$, $\forall i=1,\dots,N$, can be mathematically represented as
\begin{equation}
\label{preferenceofagents}
U_i(a,\lambda;\theta_i): \mathcal{A}  \times \Lambda \rightarrow \mathbb{R},
\end{equation}
where $\theta_i\in\Theta_i$ denotes the private information of resource agent $i$. We refer to $\theta_i$ as {\em type}. The type of resource agent $i$ may be unknown to other agents. Let $\mathcal{A}_0$ be the set of global constraints, and $\mathcal{A}_s=\mathcal{A}\times\mathcal{A}_0$ be the set of both local and global constraints on the energy allocations of the system. Then the payoff function of the coordinator agent can be mathematically represented as
\begin{equation}
\label{preferenceofcoordinator}
U_0(a,\lambda;\theta): \mathcal{A}_s  \times \Lambda \rightarrow \mathbb{R},
\end{equation}
where $\theta=(\theta_1,\ldots,\theta_N)$ and $\theta\in\Theta$ with $\Theta=\Theta_1\times\cdots\times\Theta_N$. The objectives of the resource agents and the coordinator agent are to maximize their own payoff functions, respectively.

\begin{example}
\label{exampleforsystem}
Consider the motivating example of supplying the electricity to a group of $N$ distribution loads in the presence of constrained transmission capacity. Let $\mathcal{A}_i=\{a_i|0\leq a_i^\textrm{min}\leq a_i\leq a_i^\textrm{max}\}$ denote the feasible range for the energy consumption of resource agent $i$, $\forall i=1,\dots,N$, where $a_i^\textrm{min}$ and $a_i^\textrm{max}$ are the minimally and maximally allowed energy consumption, respectively. We assume that the payoff functions of individual loads only depend on their own energy allocation and energy price. We assume a uniform price $\bar{\lambda}$ is considered across all the distribution loads, i.e., $\lambda_1=\cdots=\lambda_N=\bar{\lambda}$. Then, the payoff function of resource agent $i$, $\forall i=1,\dots,N$, is defined as
\begin{equation}
\label{agentpayoffexample}
U_i(a,\lambda;\theta_i)=U_i(a_i,\bar{\lambda};\theta_i),
\end{equation}
where $V_i(\cdot)$ is the utility of energy consumption, $a_i\in\mathcal{A}_i$ and $\bar{\lambda}\in\mathbb{R}$. In the rest of this paper, we denote the energy price as $\bar{\lambda}$ whenever it is uniform across individual resource agents.

On the other hand, the coordinator would like to supply the electricity to end users at the minimum cost while respecting the transmission capacity limit. Let $\mathcal{A}_0=\{a_i|\sum_{i=1}^Na_i\leq D\}$
denote the global constraint, where $D$ represents the energy limit for each market period due to the transmission capacity limit. Then the payoff function of the coordinator agent can be defined as
\begin{equation}
\label{coordinatorpayoffexample}
U_0(a,\lambda;\theta)=\sum_{i=1}^N U_i(a_i,\bar{\lambda};\theta_i)+\bar{\lambda}\sum_{i=1}^N a_i-C\left(\sum_{i=1}^N a_i\right),
\end{equation}
where $C(\cdot)$ is the cost for the coordinator to procure energy, $a\in\mathcal{A}_s$ and $\bar{\lambda}\in\mathbb{R}$.
\end{example}

\subsection{Control Decision}
In order to achieve maximal payoff, individual agents make control decisions to maximize their payoff functions as defined by (\ref{preferenceofagents}) and (\ref{preferenceofcoordinator}), respectively. While the payoff functions depend on the energy allocations and energy prices of all DERs, it is important to note that the control decisions of individual agents are not necessarily the energy allocations or energy prices. For instance, in some works, the coordinator determines a mapping from the energy allocation to the price \cite{li2016reverse}, \cite{ratliff2014incentive}. In some other works the agent submits a bid and the coordinator determines the energy allocation and the price according to the bid \cite{samadi2012advanced}. Therefore, the control decisions of individual agents can vary in different problems. To synthesize a unifying framework that captures these different cases, it is important  to distinguish control decisions from energy allocation and price.


Let $\gamma_i\in \Gamma_i$ be the control decision of agent $i$, $\forall i=0,1,\dots,N$, where $\Gamma_i$ is the set of constraints on control $\gamma_i$. To avoid triviality, we assume that the collective control decisions  uniquely determines the market outcome $(a,\lambda)$. With slight abuse of notation, the payoff functions (\ref{preferenceofagents}) can be written as the following function of $\{\gamma_0,\gamma_1,\dots,\gamma_N\}$,
\begin{equation}
\label{payofffunctionondecision}
U_i(\gamma_0,\gamma_1,\ldots,\gamma_N;\theta_i): \Gamma_0\times\Gamma_1\times\cdots\times\Gamma_N  \rightarrow \mathbb{R}.
\end{equation}

\begin{remark}
There exist both local and global constraints on the energy allocations for individual DERs. The local constraints, i.e., $a_i\in \mathcal{A}_i$, represent the limits at the resource level. These constraints are usually imposed by the local operational requirements. For example, the energy consumption of the air conditioner should maintain the indoor air temperature within a specified comfort range. The global constraints, i.e., $a\in \mathcal{A}_0$ represent the limits at the system level. For example, the energy allocations for individual DERs should respect the power flow constraints imposed by the network topology. In (\ref{payofffunctionondecision}), $\Gamma_i$, $\forall i=1,\ldots,N$, only represents  local constraints on $\gamma_i$. 
Note that when there are coupling constraints on control decisions in (\ref{payofffunctionondecision}), we can easily modify the proposed framework by considering the generalized equilibrium concept \cite{arrow1954existence}, \cite{debreu1952social}. This is  outside the scope of this survey. 
\end{remark}

\begin{example}
\label{exampleofcontroldecision}
Consider an example where a coordinator purchases the load shedding services from a group of responsive loads to match an electricity deficit.
To make up for the supply deficit $d$, a natural way is that the coordinator sets the load shedding price and the consumers choose the amount of load reduction according to the price. In this example, we present an alternative approach, generally referred to as the supply function bidding \cite{li2015demand}, under which the the control decisions of consumers are not the energy allocations. In the supply function bidding scheme, the coordinator requests each resource agent to submit a supply function parameterized by  $b_i$:
\begin{equation}
\label{supplyfucniton}
a_i=b_i\bar{\lambda}.
\end{equation}
A supply function maps a price to the amount of demand the resource agents are willing to reduce, and we impose that such mapping is linear with respect to $\bar{\lambda}$. In this scheme, each resource agent determines the parameter $b_i$ and is committed to the supply function for delivering $b_i\bar{\lambda}$ units of load shedding at price $\bar{\lambda}$. The price is set by the coordinator so that the total load shedding matches the electricity deficit, i.e., $\sum_{i=1}^N b_i\bar{\lambda}=d$. This leads to the following market outcome,
\begin{align}
\label{marketcleairng}
\begin{cases}
\bar{\lambda}=\dfrac{d}{\sum_{i=1}^N b_i} \\
a_i=\dfrac{b_id}{\sum_{i=1}^N b_i}, \quad \forall i=1,\ldots, N.
\end{cases}
\end{align}
In this example, the control decision of each resource agent is $b_i$. The control decision of the coordinator agent is a function that maps $(b_1,\ldots,b_N)$ to $(a,\bar{\lambda})$ as defined by (\ref{marketcleairng}). Note that the global constraint defined as $\sum_{i=1}^Na_i=d$ is always satisfied, even though resource agents only choose their bids according to local constraints, i.e., $b_i\geq 0$.
\end{example}

\subsection{Information Structure}
Information structure specifies what information is available to the agent when it makes the control decision. In transactive energy systems, information can be categorized in two classes including information on types, and information on control decisions. Information on types specifies whether one agent knows the private information of others, while information on control decision specifies whether one agent knows the control decisions of others.  Let $\mathcal{I}^{\theta}_i$ and $\mathcal{I}_i^{\gamma}$ be the sets of information available to agent $i$ on type and control decisions, respectively. These sets are determined based on the underlying assumptions of the transactive energy system. For instance, if the agent $i$ knows its own type $\theta_i$ and a joint distribution  on the types of others, which are denoted by $\Theta_{-i}(\cdot)$, then we have $\mathcal{I}_i^{\theta}=\{ \theta_i, \Theta_{-i}(\cdot) \}$. If it makes decisions after observing everyone else, then it knows the control decisions of others, thus we have $\mathcal{I}_i^{\gamma}=\{\gamma_0,\ldots,\gamma_{i-1},\gamma_{i+1},\ldots,\gamma_N \}$.

An alternative way of representing the information structure is the type dependence graph $G_{\theta}$ and the decision dependence graph $G_{\gamma}$. These graphs represent the availability of information regarding the types and control decisions of others.  The procedure of constructing $G_{\theta}$ is as follows
\begin{itemize}
\item[(1)]draw $N+1$ nodes, each of which corresponds to an agent；
\item[(2)]if agent $i$ knows $\theta_j$, draw a directed edge from $j$ to $i$；
\item[(3)]two nodes are disconnected if neither agent knows the type of the other.
\end{itemize}

\begin{figure*}[t]%
\centering
\subfigure[]{%
\label{privateinfographexample1}%
\includegraphics[width = 0.3\linewidth]{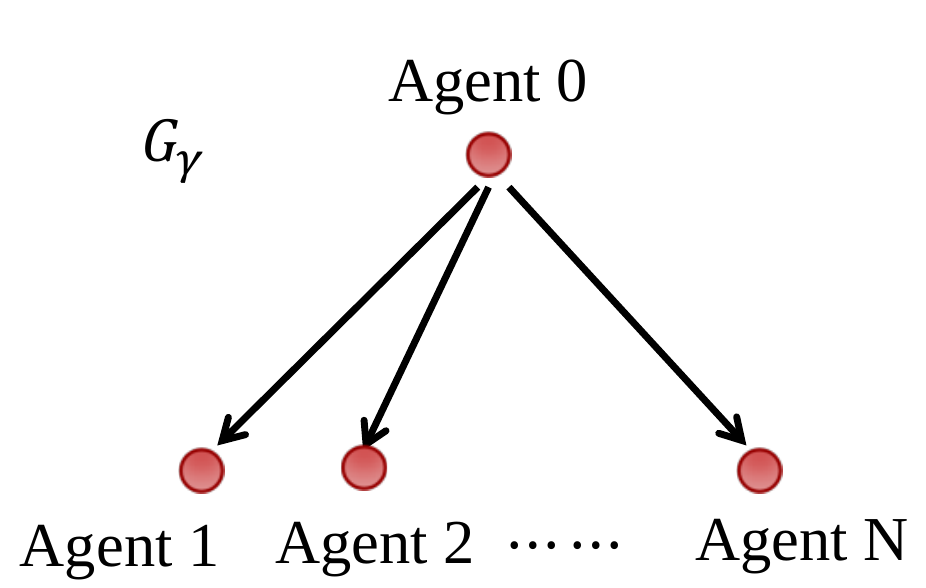}}%
\qquad
\subfigure[]{%
\label{decisiongraphexample1}%
\includegraphics[width = 0.3\linewidth]{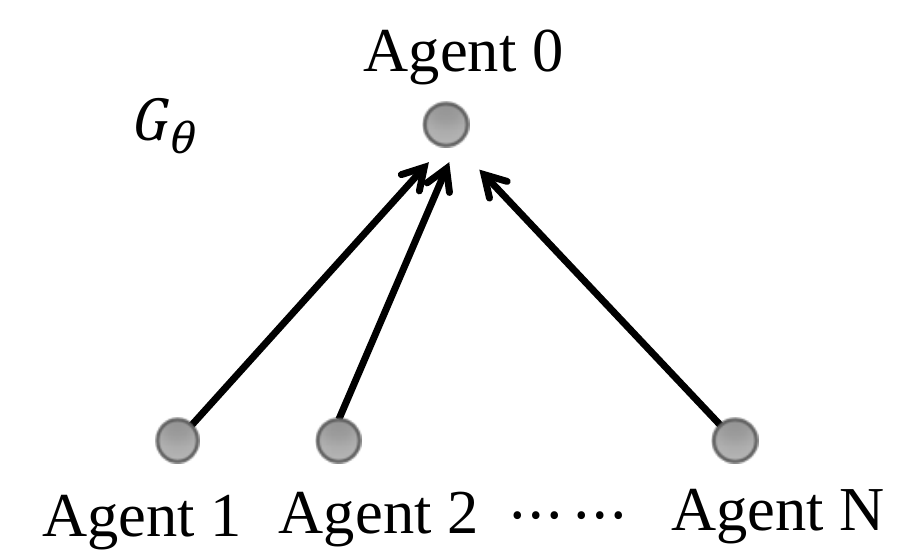}}%
\qquad
\subfigure[]{%
\label{decisiongraphexample2}%
\includegraphics[width = 0.3\linewidth]{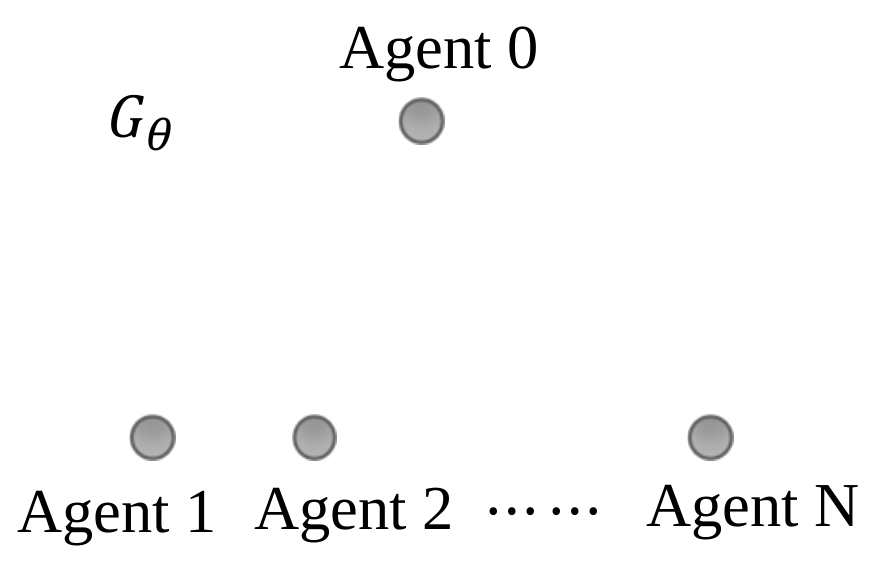}}%
\caption{The type-dependence graph and decision-dependence graph for a transactive energy system with N agent and one coordinator. (a) the coordinator knows the types of all agents, (b) the agent private information is unknown to the coordinator, (c) the coordinator makes control decision first, followed by agents who observe the coordinator decision and act subsequently. }
\end{figure*}

A graph constructed by this procedure captures the private information structure of the transactive energy system. As an example, we show two type dependence graphs in Figures \ref{privateinfographexample1} and \ref{decisiongraphexample1}. In Figure \ref{privateinfographexample1}, the coordinator knows the type of all agents, but agents do not know the types of other agents.  In Figure \ref{decisiongraphexample1}, the coordinator does not know the type of any agent. Therefore, there is no edge in the graph.

The decision-dependence graph is closely related to the order of decision making: If all agents made decisions simultaneously, then no one knows the decisions of others. In our paper, we group the agents in a number of stages, and pass the decisions of agents in earlier stages to all agents in the subsequent stages. In this case, all agents know the decision of the agents taking actions before them.
 The detailed procedure for constructing $G_{\gamma}$ is  as follows:
\begin{itemize}
\item[(1)]divide all agents into $M$ groups, denoted as $S_1,\ldots,S_M$.
\item[(2)]if $i\in S_m$, then agent $i$ makes control decision at stage $m$.
\item[(3)]$\forall i\in S_{m}, \forall j\in S_{p}$, agent $i$ knows $\gamma_j$ if and only if $m>p$.
\item[(4)]if agent $i$ knows $\gamma_j$ before making a control decision, draw a directed edge from $j$ to $i$.
\end{itemize}
Following the above procedure, a decision dependence graph can be obtained to capture the information structure for sequential decision making. As an example, Figure \ref{decisiongraphexample2} shows a transactive energy system where the coordinator determines $\gamma_0$ in the first stage, and then passes it to  agents making decisions in the second stage.


\subsection{Solution Concept}
The payoff function  of each agent  (\ref{payofffunctionondecision}) depends on the control decisions of other agents. Therefore, the solution to the transactive energy system is an equilibrium. An equilibrium point predicts the stable outcomes of the decision making process: at the equilibrium, no agent is willing to deviate given the information available to it.  The solution concept for the  equilibrium can be complicated due to the sequential nature of  decision making.
For ease of presentation, we first introduce the solution concepts for a singe-stage problem, then we extend these solution concepts to  multi-stage problems.

\subsubsection{Single-Stage Transactive Energy Systems}
In a single-stage transactive energy system, all the agents make the control decisions simultaneously. If their payoff functions are decoupled, it becomes a standard optimization. Otherwise, it is a game, where the most widely used solution concept is the Nash equilibrium.
\begin{definition}
\label{Nashdefinition}
A strategy profile $(\gamma_0^*,\gamma_1^*,\ldots,\gamma_N^*)$ is a Nash equilibrium if for all $i=0,\ldots,N$ and all $\gamma_i\in \Gamma_i$, we have
\begin{equation}
U_i(\gamma_i^*,\gamma_{-i}^*; \theta_i)\geq  U_i(\gamma_i,\gamma_{-i}^*; \theta_i),
\end{equation}
where $\gamma_{-i}^*$ represents the decisions of all agents except sagent $i$, i.e., $\gamma_{-i}^*=\{\gamma_0^*,\ldots,\gamma_{i-1}^*,\gamma_{i+1}^*,\ldots,\gamma_N^*\}$.
\end{definition}
At Nash equilibrium, it is to the agent's best benefit to follow the Nash equilibrium strategy given that all other agents also follow such equilibrium strategy. Therefore, no one will unilaterally deviate if the system has arrived at the Nash equilibrium.

A relaxed version of Nash equilibrium is often used in large population games \cite{huang2007large}, \cite{li2017connections}. It is referred to as the $\epsilon$-Nash equilibrium.
\begin{definition}
\label{relaxNashdefinition}
A strategy profile $(\gamma_0^*,\gamma_1^*,\ldots,\gamma_N^*)$ is an $\epsilon$-Nash equilibrium if for all $i=0,\ldots,N$ and all $\gamma_i \in \Gamma_i$,  we have
\begin{equation}
U_i(\gamma_i^*,\gamma_{-i}^*; \theta_i)\geq  U_i(\gamma_i,\gamma_{-i}^*; \theta_i)-\epsilon,
\end{equation}
where $\gamma_{-i}^*=\{\gamma_0^*,\ldots,\gamma_{i-1}^*,\gamma_{i+1}^*,\ldots,\gamma_N^*\}$.
\end{definition}

The $\epsilon$-Nash equilibrium implicitly assumes that each agent is indifferent with respect to a small change of $\epsilon$ in its payoff. Under this assumption, an agent is motivated to play an $\epsilon$-Nash equilibrium if other agents also play an $\epsilon$-Nash equilibrium. If $\epsilon=0$, then $\epsilon$-Nash equilibrium strategy  becomes an exact Nash equilibrium.

At a Nash equilibrium, agents will not deviate if others are also at the Nash equilibrium. However, this concept is hard to realize when agents have private information. In this case,  each agent does not know the type of others and can not predict the actions of others in order to make decisions on its own. Addressing this issue requires a stronger solution concept, such as the Bayesian Nash equilibrium and the dominant strategy equilibrium.

If all agents have a prior distribution on the types of others, they can predict the behaviors of other agents based on the prior distribution and make control decision accordingly. This can be captured by the Bayesian Nash equilibrium. To this end,  we explicitly write $\gamma_i$ as $\gamma_i(\theta_i)$ to emphasize that the control decision depends on the private information. Formally, a Bayesian Nash equilibrium is defined as follows.

\begin{definition}
\label{BayesianNashdefinition}
A strategy profile $(\gamma_0^*(\cdot),\gamma_1^*(\cdot),\ldots,\gamma_N^*(\cdot))$ is a Bayesian Nash equilibrium if for all $i=0,\ldots,N$ and all $\gamma_i(\cdot)\in \Gamma_i$, we have
\begin{equation}
\mathbb{E}_{\theta_{-i}} \left[ U_i(\gamma_i^*(\cdot),\gamma^*_{-i}(\cdot); \theta_i)|\theta_i\right]\geq  \mathbb{E}_{\theta_{-i}}\left[ U_i(\gamma_i(\cdot),\gamma^*_{-i}(\cdot); \theta_i)|\theta_i\right],
\end{equation}
\end{definition}
where $\theta_{-i}$ is the types of all agents except agent $i$. The expectation is over the prior distribution on $\theta_{-i}$. 

When the prior distribution on types is not available,  the agents have to make decisions without using any information about others. This can be captured by a stronger solution concept, defined as follows \cite{diamantaras2009toolbox}, \cite{mas1995microeconomic}.
\begin{definition}
\label{dominantdefinition}
A strategy profile $(\gamma_0^*,\gamma_1^*,\ldots,\gamma_N^*)$ is a dominant strategy equilibrium if for all $i=0,\ldots,N$ and all $\gamma_i\in \Gamma_i$,  we have
\begin{equation}
U_i(\gamma_i^*,\gamma_{-i}; \theta_i)\geq  U_i(\gamma_i,\gamma_{-i}; \theta_i),
\end{equation}
for all feasible $\gamma_{-i}$, where $\gamma_{-i}=\{\gamma_0,\ldots,\gamma_{i-1},\gamma_{i+1},\ldots,\gamma_N\}$.
\end{definition}

At a dominant strategy equilibrium, a rational agent always follows the equilibrium strategy regardless the decisions of others. In contrast,  at Nash equilibrium or Bayesian Nash equilibrium, the agent plays the equilibrium strategy only when it has a correct forecast of the actions of other agents. Therefore, compared to (Bayesian) Nash equilibrium, the dominant strategy equilibrium is much easier to implement: each agent simply chooses the dominant strategy equilibrium without knowing anything about others. However, it is also much harder to guarantee  the existence of  dominant strategy equilibrium.

\subsubsection{Multi-Stage Transactive Energy Systems}
The solution concept for a multi-stage transactive energy system is more complicated and carries heavy notation \cite{bacsar1998dynamic}. For ease of presentation, we start with a two-stage problem with control-dependence graph depicted in Figure \ref{decisiongraphexample2}. For this case, the coordinator agent determines $\gamma_0$ at the first stage, then each resource agent chooses $\gamma_i$ according to $\gamma_0$ in the second stage. Given the first stage decision $\gamma_0$, the second-stage is reduced to a single-stage problem. Therefore, a two-stage problem is equivalent to single-stage problem parameterized by the first-stage decision. Under a proper solution concept for the second stage, we denoted the solution as $\gamma_i^*(\gamma_0)$. Note that $\gamma_i^*$ depends on  $\gamma_0$ as different $\gamma_0$ leads to different games and corresponding equilibria.

Given the response of resource agents in the second stage, the first-stage decision making of the coordinator agent can be formulated as the following optimization problem,
\begin{IEEEeqnarray}{r'l}
\max_{\gamma_0\in \Gamma_0} & U_0(\gamma_0,\gamma_1^*(\gamma_0),\ldots,\gamma_N^*(\gamma_0)). \IEEEyesnumber \label{stage0solution}
\end{IEEEeqnarray}

In this two-stage game, the solution concept of the second stage is a game equilibrium, while the solution concept of the first stage is just a standard optimization. This is because no one competes with the coordinator in the first stage. If there are multiple coordinators, both stages may have the game equilibria as the solution concept \cite{maharjan2013dependable}.

\begin{remark}
The solution concept can be defined similarly for games with more than two stages. In this case, the decision of each agent is   parameterized by the decisions of all agents making decisions prior to it, and the corresponding solution concept can be chosen for each stage of the game.
\end{remark}



\begin{figure*}[t]
\centering
\includegraphics[width = 0.75\linewidth]{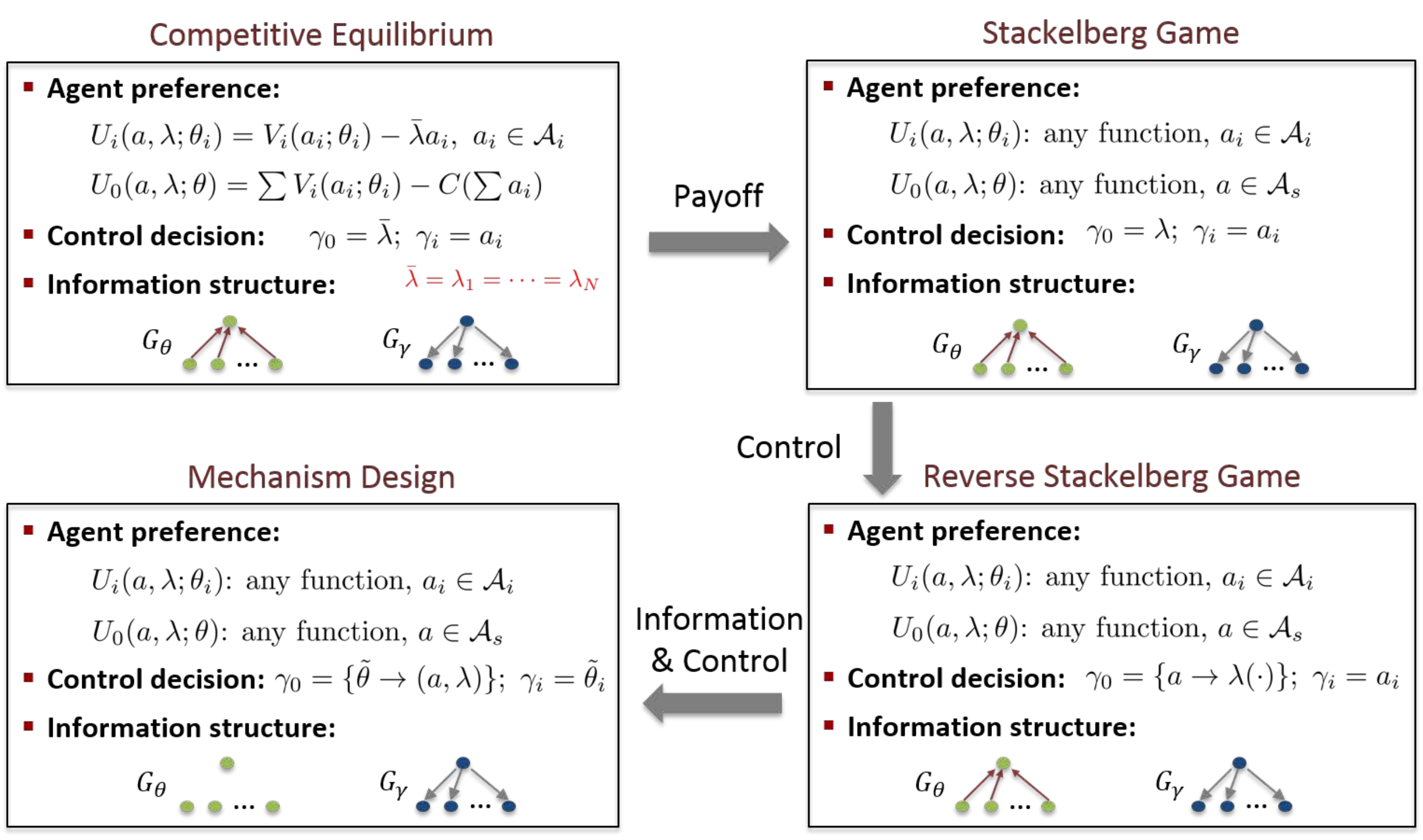}
\caption{Comparisons for different categories of transactive energy systems. Competitive equilibrium is different from Stackelberg games in coordinator's payoff function, Stackelberg game differs from Reverse Stackelberg game in control decisions, and mechanism design differs from reverse Stackelberg game in both information structure and control decisions. }
\label{diagram_comparison}
\end{figure*}

In the rest of this paper, we present several representative transactive energy systems for solving the same market-based coordination problem with different assumptions, aiming to clearly illustrate the involved subtleness and complexity. Each of these transactive energy systems has its own special formulation under the proposed unifying framework, serving as a typical example of a class of transactive energy systems that have been extensively studied in the literature. The comparisons between these representative systems are discussed with detailed literature survey and a pictorial summary in Figure~\ref{diagram_comparison}.

\section{Competitive Equilibrium}
\label{CEsection}
Efficient energy allocation is the main objective of many works \cite{nguyen2012walrasian}, \cite{papadaskalopoulos2013decentralized, hansen2015dynamic, chen2012optimal, li2011optimal, jokic2010price, singh2016decentralized}. Although these works differ in various aspects, they share some common features related to the well-established concept of competitive equilibrium in microeconomics \cite{mas1995microeconomic}.
In the following, we synthesize a transactive energy system as a typical example for this class of works.

\subsection{Transactive Energy System}
Consider the problem of market-based coordination for a group of $N$ distribution loads. Note that, for the simplicity of the following discussions, only distribution loads are considered. However, the results presented herein hold when other types of DERs are considered. In this section, we assume that the payoff functions of individual loads depend on their own energy allocation and energy price. Furthermore, their payoff functions are quasi-linear with respect to the energy price. We also assume that the energy price is uniform for all loads. Under these assumptions, the objective of resource agent $i$, $\forall i=1,\dots,N$, is to maximize the following payoff function,
\begin{equation}
\label{CEagents}
U_i(a,\lambda;\theta_i)=V_i(a_i;\theta_i)-\bar{\lambda} a_i,
\end{equation}
where $V_i(\cdot)$ is the utility of energy consumption, $a_i\in\mathcal{A}_i$, and $\bar{\lambda}\in\mathbb{R}$. The objective of the coordinator is to obtain efficiency energy allocations that maximize the following payoff function,
\begin{equation}
\label{CEcoordinatorpayoof}
U_0(a,\lambda;\theta)=\sum_{i=1}^N V_i(a_i;\theta_i)-C\left(\sum_{i=1}^N a_i\right),
\end{equation}
where $C(\cdot)$ is the cost for the coordinator to purchase energy from outside the system, and $a\in\mathcal{A}$.

We assume that the coordinator agent may or may not know the private information of individual resource agents. It determines the energy price $\bar{\lambda}$ in the first stage, and then resource agents choose their energy consumption $a_i$ in the second stage according to $\bar{\lambda}$. Hence, this is a transactive energy system with the following formulation,
\begin{itemize}
\item agent preference: the resource agents and the coordinator agent have preferences defined by (\ref{CEagents}) and  (\ref{CEcoordinatorpayoof}), respectively.
\item control decision: $\gamma_0=\bar{\lambda}$ and $\gamma_i=a_i, \forall i=1,\ldots,N$.
\item information structure: the type-dependence graph as in Figure \ref{privateinfographexample1} or \ref{decisiongraphexample1}, and decision-dependence graph as in Figure \ref{decisiongraphexample2}.
\item solution concept: given $\bar{\lambda}$, the payoff functions of resource agents are decoupled, so it is a standard optimization problem.
\end{itemize}
It follows from (\ref{CEagents}) and  (\ref{CEcoordinatorpayoof}) that the optimal energy price $\bar{\lambda}^*$ of the coordinator agent is obtained by solving the following optimization problem,
\begin{IEEEeqnarray}{r'l}
\label{transactiveenergysystem}
\max_{\bar{\lambda}\in\mathbb{R}} & \sum_{i=1}^N V_i\left(a_i^*(\bar{\lambda});\theta_i\right)-C\left(\sum_{i=1}^N a_i^*(\bar{\lambda})\right) \IEEEyesnumber \IEEEyessubnumber \label{transactiveenergysystema}\\
\text{s.t} & a_i^*(\bar{\lambda})=\argmax_{a_i\in\mathcal{A}_i} V_i(a_i;\theta_i)-\bar{\lambda} a_i.
\IEEEyessubnumber \label{transactiveenergysystemb}
\end{IEEEeqnarray}

This transactive energy system involves sequential decision making: the coordinator determines the optimal energy price first, then individual loads choose their optimal energy consumption accordingly. In this case, the coordinator should anticipate the responses of individual loads to its decision, and exploit this knowledge to determine the price. Such problem is often referred to as bi-level problem \cite{gabriel2012complementarity}, where the coordinator's problem has another optimization problem as a constraint \cite{vicente1994bilevel}, \cite{dempe2003annotated}, \cite{colson2007overview}.

\subsection{Connection to Competitive Equilibrium}

In this subsection, the competitive equilibrium \cite{mas1995microeconomic} is first introduced. Then, its connection to the transactive energy system discussed in the last subsection is presented. We show that the solution to  (\ref{transactiveenergysystem}) is essentially a competitive equilibrium.

Consider a competitive market  with  one supplier and $N$ consumers trading electricity. Individual consumers want to maximize the following payoff functions,
\[
U_i(a_i,\lambda;\theta_i)=V_i(a_i;\theta_i)-\bar{\lambda} a_i, \quad \forall i=1,\dots,N,
\]
where $a_i\in\mathcal{A}_i$ and $\bar{\lambda}\in\mathbb{R}$, while the supplier wants to maximize the following payoff function,
\[
\label{eq:CEpreferencesupplier}
U_s(y,\lambda;\theta_s)=\bar{\lambda} y-C(y;\theta_s),
\]
where $\theta_s$ denotes the private information of the supplier, $y\in\mathcal{Y}$ is the energy generation, and $\mathcal{Y}$ is the feasible range for $y$. In a competitive market, both the consumer and the supplier take $\bar{\lambda}$ as given, and they choose consumption or production level to maximize their payoffs. The market outcome is the competitive equilibrium, formally defined as follows \cite[Chap 10]{mas1995microeconomic}:
\begin{definition}
\label{definitionofCE}
The energy allocation vector $(a_1^*,\ldots,a_N^*,y^*)$ and the energy price $\bar{\lambda}^*$ constitute a competitive equilibrium if the following conditions are satisfied.
\begin{itemize}
\item[(i)] Each consumer maximizes the payoff at the equilibrium, i.e., $a_i^*$ is the optimal solution to the following problem,
\begin{IEEEeqnarray}{r'l}
\max_{a_i} & V_i(a_i;\theta_i)-\bar{\lambda}^* a_i \IEEEyesnumber \IEEEyessubnumber \nonumber\\
\text{s.t} & a_i\in \mathcal{A}_i. \IEEEyessubnumber \nonumber
\end{IEEEeqnarray}

\item[(ii)] The supplier maximizes the payoff at the equilibrium, i.e., $y^*$ is the optimal solution to the following problem,
\begin{IEEEeqnarray}{r'l}
\max_{y} & \bar{\lambda}^* y-C(y;\theta_s) \IEEEyesnumber \IEEEyessubnumber \nonumber\\
\text{s.t} & y\in \mathcal{Y}. \IEEEyessubnumber \nonumber
\end{IEEEeqnarray}

\item[(iii)] Total supply and total demand are balanced, i.e., $\sum_{i=1}^N a_i^*=y^*$.
\end{itemize}
\end{definition}

In the above, conditions (i) and (ii) state that individual agents must take the energy price as given and choose the energy consumption (or production) to maximize their payoff functions. Condition (iii) requires that at the equilibrium point, the overall consumption and production identified by (i) and (ii) must be mutually compatible. That is, no excess demand or supply appears at the equilibrium point.
It is well-known that competitive equilibrium maximizes the social welfare \cite[Chap 10]{mas1995microeconomic},
\begin{proposition}
\label{CEmaxsocialeflare}
The competitive equilibrium $(a_1^*,\ldots,a_N^*,y^*,\bar{\lambda}^*)$ is the optimal solution to the following social welfare optimization problem,
\begin{IEEEeqnarray}{r'l}
\label{socialwelfareopt}
\max_{a_1,\ldots,a_N,y} & \sum_{i=1}^N V_i(a_i;\theta_i)-C(y;\theta_s) \IEEEyesnumber \IEEEyessubnumber \label{socialwelfareopta}\\
\text{s.t} & \sum_{i=1}^N a_i =y \IEEEyessubnumber \label{socialwelfareoptb} \\
& y\in \mathcal{Y} \IEEEyessubnumber \label{socialwelfareoptc} \\
& a_i\in \mathcal{A}_i,\ \forall i=1,\ldots,N. \IEEEyessubnumber \label{socialwelfareoptd}
\end{IEEEeqnarray}

\end{proposition}

The proof of Proposition (\ref{CEmaxsocialeflare}) can be found in \cite[Chap 10]{mas1995microeconomic}. We can show that the competitive equilirbium is the solution to the transactive energy system defined by (\ref{CEagents})--(\ref{transactiveenergysystem}).


\begin{proposition}
\label{CEtheorem}
Assume that $V_i(\cdot; \theta_i)$ is concave, $C(\cdot)$ is convex. Assume that $\mathcal{A}_i$ and $\mathcal{Y}$ are convex and have at least one interior point. Then   $(a_1^*,\ldots,a_N^*,y^*,\bar{\lambda}^*)$ is a competitive equilibrium if and only if $(a_1^*,\ldots,a_N^*,\bar{\lambda}^*)$ is the optimal solution to (\ref{transactiveenergysystem}).
\end{proposition}

Proposition \ref{CEtheorem} implies that the optimal solution to  (\ref{transactiveenergysystem}) is a competitive equilibrium. It can be proved by viewing $\bar{\lambda}^*$ as the Lagrange multiplier of (\ref{socialwelfareoptb}). Together with the results of Proposition \ref{CEmaxsocialeflare}, we have the following corollary.

\begin{corollary}
\label{CEequivalencetoSWM}
Under the assumptions of Proposition  \ref{CEtheorem}, the optimal solution to (\ref{transactiveenergysystem}) is also the optimal solution to the social welfare optimization problem (\ref{socialwelfareopt}).
\end{corollary}


\begin{remark}
The concept of competitive equilibrium has been well-studied. The social welfare optimization problem (\ref{socialwelfareopt}) can be efficiently solved in many different ways as discussed in the following subsection. The main purpose of this subsection is to establish the connection of the transactive energy system (\ref{CEagents})--(\ref{transactiveenergysystem}) to the well-established concept of competitive equilibrium so that all the results on competitive equilibrium can be directly applied to a large class of transactive energy problems \cite{chen2010two}, \cite{nguyen2012walrasian}, \cite{papadaskalopoulos2013decentralized}, \cite{hansen2015dynamic}, \cite{chen2012optimal}, \cite{li2011optimal}, \cite{jokic2010price}, \cite{singh2016decentralized}, \cite{yang2017distributed}, \cite{dominguez2012decentralized},  \cite{zhang2012convergence}.
\end{remark}

\subsection{Computation of Competitive Equilibrium}
The social welfare optimization problem (\ref{socialwelfareopt}) is a simple convex optimization problem if the coordinator agent knows the types of all the resource agents. However, if the coordinator does not know these private information, the objective function of (\ref{socialwelfareopt}) is unknown. In the literature, both non-iterative and  iterative approaches have been proposed to address this issue.

The non-iterative approach adopts an auction-based market to determine the optimal energy price, referred to as market clearing, in a hierarchical manner \cite{Fuller11}, \cite{Hammerstrom07}, \cite{Widergren14Report}, \cite{Widergren14ISGT}, \cite{li2016market}, \cite{Lian2018metrics}. The details of the auction-based approach is summarized in Algorithm~\ref{hierarchicalalgorithm}. To implement this approach, individual resource agents have to determine their demand curves as bids at the beginning of each market period and send them to the coordinator agent. After receiving all the demand bids, the coordinator aggregates all the demand curves, and find the intersection between the supply curve and the aggregated demand curve.

\begin{algorithm}
\caption{The Auction-based Approach to Compute the Competitive Equilibrium} \label{hierarchicalalgorithm}
\begin{algorithmic}[1]

\REQUIRE Transactive Energy System (\ref{CEagents}) and (\ref{CEcoordinatorpayoof}). 

\STATE Individual resource agents determine their demand curves by solving
\begin{equation}
\label{hierarchical1}
a_i(\bar{\lambda})=\argmax_{a_i\in \mathcal{A}_i} V_i(a_i;\theta_i)-\bar{\lambda} a_i.
\end{equation}

\STATE Individual resource agents send their demand curves $a_i(\bar{\lambda})$ as demand bids to the coordinator agent.

\STATE The Coordinator agent computes the supply curve by solving
\begin{equation}
\label{hierarchical2}
y(\bar{\lambda})=\argmax_{y\in \mathcal{Y}} \bar{\lambda}y-C(y;\theta_s).  
\end{equation}

\STATE The coordinator agent determine the energy price so that
\begin{equation}
\label{hierarchical3} 
\sum_{i=1}^N a_i(\bar{\lambda})=y(\bar{\lambda}).
\end{equation}

\ENSURE the solution $(a_1^*,\ldots,a_N^*,\bar{\lambda}^*)$.
\end{algorithmic}
\end{algorithm}

The iterative approach determines the market clearing price in an iterative manner. Many different algorithms have been proposed in the literature, including ratio consensus \cite{yang2017distributed}, \cite{dominguez2012decentralized}, average consensus \cite{zhang2012convergence},  primal-dual algorithms  \cite{chen2010two}, \cite{hansen2015dynamic}, among others. In this subsection, we present an algorithm based on primal-dual iterations. This algorithm enables parallel implementation and scales well with the size of the transactive energy system. Therefore, it has been widely adopted in many works \cite{chen2010two}, \cite{nguyen2012walrasian}, \cite{papadaskalopoulos2013decentralized}, \cite{hansen2015dynamic}, \cite{chen2012optimal}, \cite{li2011optimal}, \cite{jokic2010price}, \cite{singh2016decentralized}. The details of the primal-dual algorithm is summarized in Algorithm~\ref{primaldualalgorithm}. It is also referred to as Walrasian auction \cite{mas1995microeconomic}, and its solution is the Walrasian equilibrium. To implement this algorithm, the coordinator starts with an initial guess for the energy price $\bar{\lambda}$ and broadcasts it to all the resource agents. Each resource agent then solves the individual payoff maximization problem according to (\ref{primaldual1}). The solutions to the payoff maximization problems are used to update the dual $\bar{\lambda}$ according to (\ref{primaldual3}). The updated dual variable is then broadcast to resource agents again and this procedure is iterated until the energy price converges. This iterative algorithm is a sub-gradient algorithm that converges to the optimal energy allocation and energy price pair under the assumptions of Proposition \ref{CEtheorem} \cite{boyd2011distributed}. However, we would like to point out that these algorithms are not strategy-proof in that the resource agents may not necessarily have incentives to follow the specified iteration protocol.

\begin{algorithm}
\caption{The Primal-Dual Algorithm to Compute the Competitive Equilibrium} \label{primaldualalgorithm}
\begin{algorithmic}[1]

\REQUIRE Transactive Energy System (\ref{CEagents}) and (\ref{CEcoordinatorpayoof}).

\STATE The coordinator agent generates an initial guess of the energy price, $\bar{\lambda}(0)$.

\WHILE {Stopping criteria not satisfied}

\STATE Individual resource agents determine control decisions by solving
\begin{equation}
\label{primaldual1}
a_i(k)=\argmax_{a_i\in \mathcal{A}_i} V_i(a_i;\theta_i)-\bar{\lambda}(k-1) a_i.  
\end{equation}

\STATE The Coordinator agent computes the following value,
\begin{equation}
\label{primaldual2}
y(k)=\argmax_{y\in \mathcal{Y}} \bar{\lambda}(k-1) y-C(y;\theta_s)  
\end{equation}

\STATE The coordinator agent updates the energy price,
\begin{equation}
\label{primaldual3} 
\bar{\lambda}(k) =\bar{\lambda}(k-1)+\gamma(k)  \left(\sum_{i=1}^N a_i(k)-y(k)\right).
\end{equation}

\ENDWHILE 

\ENSURE the solution $(a_1^*,\ldots,a_N^*,\bar{\lambda}^*)$.
\end{algorithmic}
\end{algorithm}  


\section{Stackelberg Game}
For the transactive energy system defined by (\ref{CEagents})--(\ref{transactiveenergysystem}), the payoff functions of individual loads are assumed to be quasi-linear with the energy price.  In this section, we generalize the structure of the payoff functions and cast the coordination problem as a Stackelberg game instead. Existing results on this class of transactive energy systems are surveyed with challenges clearly identified.

\subsection{Transactive Energy System}
Consider again the problem of market-based coordination for a group of $N$ distribution loads. In this section, however, we assume the following general payoff functions without imposing any specific forms as in the previous section for individual loads and the coordinator, respectively, 
\begin{equation}
\label{agentpayoffstackel}
U_i(a,\lambda;\theta_i), \quad \forall i=1,\dots,N,
\end{equation}
where $a_i\in \mathcal{A}_i$ and $\lambda\in\Lambda$, and
\begin{equation}
\label{coordintorpayoffstackel}
U_0(a,\lambda;\theta),
\end{equation}
where $a\in \mathcal{A}_s$ and $\lambda\in\Lambda$. We assume that the coordinator knows all the private information. The coordinator agent announces a price vector $\lambda\in \Lambda$ in the first stage, then each resource agent chooses $a_i$ accordingly in the second stage. In determining the energy price, the coordinator should anticipate the resource responses to its price so as to maximize its payoff (\ref{coordintorpayoffstackel}). In the meanwhile, each resource agent also chooses $a_i$ to solve (\ref{agentpayoffstackel}). This class of problems are widely studied in the literature \cite{maharjan2013dependable},  \cite{tushar2012economics}, \cite{asimakopoulou2013leader}, \cite{zhong2013coupon}, \cite{maharjan2016demand},  and we can summarize it as the following transactive energy system,
\begin{itemize}
\item agent preference: the resource agents and the coordinator agent have preferences defined by (\ref{agentpayoffstackel}) and (\ref{coordintorpayoffstackel}), respectively.
\item control decision: $\gamma_0=\lambda$ and $\gamma_i=a_i$ for $\forall i=1,\ldots,N$.
\item information structure: the type-dependence graph  as in Figure \ref{privateinfographexample1}, and decision-dependence graph as in Figure \ref{decisiongraphexample2}.
\item solution concept: given $\lambda$, the decisions of resource agents are coupled, and the solution concept is the Nash equilibrium.
\end{itemize}

Given the energy price $\lambda$, the problem of the second stage is a game due to the coupling among payoff functions of resource agents. We define the game equilibrium as
\begin{equation*}
\label{lowerlvelstackel}
a_i^*(\lambda)=\argmax_{a_i\in A_i} U_i(a_i,a_{-i}^*,\lambda;\theta_i), \quad \forall i=1,\ldots,N.
\end{equation*}
where $a_{-i}^*=(a_1^*,\ldots,a_{i-1}^*,a_{i+1}^*,\ldots,a_N^*)$, and $a^*(\lambda)$ is the optimal response of each resource agent to the energy price $\lambda$. Then the decision of the coordinator agent reduces to the solution to the following optimization problems,
\begin{IEEEeqnarray}{r'l}
\label{stackelberggame}
\max_{\lambda\in \Lambda} & U_0(a_1^*(\lambda),\ldots,a_N^*(\lambda), \lambda;\theta) \IEEEyesnumber \IEEEyessubnumber \label{stackelberggamea}\\
\text{s.t} & a_i^*(\lambda)=\argmax_{a_i\in A_i} U_i(a_i,a_{-i}^*,\lambda;\theta_i).  \IEEEyessubnumber \label{stackelberggameb}
\end{IEEEeqnarray}

The transactive energy system defined by (\ref{agentpayoffstackel})--(\ref{stackelberggame}) is different from the one defined by (\ref{CEagents})--(\ref{transactiveenergysystem}). In the previous section, the energy price $\lambda$ is uniform across all resource agents, that is, $\lambda_i=\bar{\lambda}$. Furthermore, the payoff functions of resource agents (\ref{CEagents}) are quasi-linear, and only depend on $a_i$ and $\bar{\lambda}$. In this section, the energy price is not necessarily uniform, and the payoff functions of resource agents (\ref{agentpayoffstackel}) are general functions of $a$ and $\lambda$. For this reason, the solution concept of the second stage changes from a standard optimization to the Nash equilibrium.

\begin{example}
\label{example3}
Consider again the motivating example of supplying the electricity to a group of $N$ distribution loads in the presence of constrained transmission capacity. The coordinator agent can indirectly affect the end users by setting retail prices. Different from Example \ref{exampleforsystem}, the structures of the payoff functions of the resource agents and the coordinator agent can be fairly general. We assume that the coordinator agent knows all the private information $\theta$. It determines the energy price $\bar{\lambda}$ in the first stage, and each resource agent determines $a_i$ accordingly in the second stage. Given the energy price $\bar{\lambda}$, the optimal responses of resource agents are given as
\begin{equation*}
\label{agentresponse}
 a_i^*(\bar{\lambda})=\argmax_{a_i\in\mathcal{A}_i} V_i(a_i;\theta_i)-\bar{\lambda}a_i, \quad \forall i=1,\ldots,N.
\end{equation*}
Taking into account the optimal responses of resource agents, the coordinator agent can determine the optimal energy price $\bar{\lambda}$ by formulating the following optimization problem,
\begin{IEEEeqnarray}{r'l}
\max_{\bar{\lambda}\in\mathbb{R}} & U_0(a^*(\bar{\lambda}),\bar{\lambda};\theta) \IEEEyesnumber \label{stackelberggameexample} \nonumber\\
\text{s.t} & \sum_{i=1}^N a_i^*(\bar{\lambda})\leq D \nonumber \\
& a_i^*(\bar{\lambda})=\argmax_{a_i\in\mathcal{A}_i} V_i(a_i;\theta_i)-\bar{\lambda}a_i, \quad \forall i=1,\ldots,N, \nonumber
\end{IEEEeqnarray}
where $a^*(\bar{\lambda})=(a_1^*(\bar{\lambda}),\ldots,a_N^*(\bar{\lambda}))$. This is a special case of (\ref{stackelberggame}).
\end{example}

\subsection{State of the Art and Challenges}
Similar to the transactive energy system defined by (\ref{CEagents})--(\ref{transactiveenergysystem}), the transactive energy systems defined by (\ref{agentpayoffstackel})--(\ref{stackelberggame}) also involves sequential decision making. Thus, it is also a bi-level problem. In the previous section, it follows from Corollary (\ref{CEequivalencetoSWM}) that the bi-level optimization problem (\ref{transactiveenergysystem}) can be transformed into the a single-level optimization, that is, social welfare optimization problem (\ref{socialwelfareopt}), which can be efficiently solved. However, this simplification is only possible for  quasi-linear payoff functions with local constraints. These properties do not carry over to more general cases. Therefore, comparing to (\ref{transactiveenergysystem}), the optimization problem (\ref{stackelberggame}) is more general, but also significantly more challenging from a computational viewpoint.

In the literature, problem (\ref{stackelberggame}) is often referred to as a Stackelberg game  \cite{gabriel2012complementarity}, \cite{von1952theory}, \cite{simaan1973stackelberg}. The concept of Stackelberg game originates from economics \cite{von1952theory}, and has received substantial attentions in recent years in the problem of coordinating DERs  \cite{maharjan2013dependable},  \cite{yu2016real},    \cite{tushar2012economics}, \cite{zhong2013coupon} ,\cite{asimakopoulou2013leader}, \cite{maharjan2016demand}, \cite{coogan2013energy}, \cite{tushar2014prioritizing},  \cite{tushar2014energy}, \cite{tushar2015three}.  Different from the standard optimization, it has an inner optimization problem as a constraint. For this reason, Stackelberg game generally involves non-convex and non-differentiable optimization problems, which is rather challenging from a computational viewpoint. In fact, problem (\ref{stackelberggame}) is strongly NP-hard \cite{hansen1992new}. Furthermore, the problem of merely checking the global or local optimality for (\ref{stackelberggame}) is NP-hard \cite{vicente1994bilevel}, even if all payoff functions and constraints are linear.

Since the non-convexity and non-differentiability of (\ref{stackelberggame}) greatly undermine our ability to solve the general problem for the global optimality, special attentions have been given to the simpler cases where all payoff functions and constraints are linear \cite{wen1991linear}, \cite{bard2013practical}. For those problems with linear payoff functions, the optimal solution can be found at the vertex of the underlying polyhedron. This leads to many global algorithms based on enumerative computation such as vertex enumeration \cite{campelo2000note}, \cite{candler1982linear}, \cite{onal1993modified}, branch-and-bound algorithms \cite{bard1982explicit}, \cite{fortuny1981representation}, complementary pivoting algorithms \cite{judice1988solution}, \cite{judice1992sequential}, among others. These algorithms have been extended to  quadratic payoff functions, where the globally optimal solution can be also guaranteed \cite{judice1994linear}, \cite{vicente1994descent}, \cite{bard1990branch}. On the other hand, if the objectives are general nonlinear functions, solution attempts include  gradient descent \cite{kolstad1990derivative}, \cite{savard1994steepest}, penalty function method \cite{shimizu1981new}, \cite{aiyoshi1981hierarchical}, \cite{aiyoshi1984solution}, trust region algorithms \cite{conn2000trust}, \cite{liu1998trust}, etc. However, due to the non-convexity of (\ref{stackelberggame}), most of these algorithms can only find locally optimal solutions.

\section{Reverse Stackelberg Game}
In many problems of interest, the control decision of the coordinator agent is not a price, but rather a pricing function that depends on the control decisions of resource agents \cite{moulin2010efficient}, \cite{barreto2013design}, \cite{barreto2014incentives}. This leads to a new class of transactive energy systems. In this section, we formulate the coordination problem as a reverse Stackelberg game. We survey the existing results on this class of transactive energy
systems and identify the associated technical challenges.

\subsection{Transactive Energy System}
Consider again the market-based coordination problem with the general payoff functions as in (\ref{agentpayoffstackel}) and (\ref{coordintorpayoffstackel}). In this section, however, we assume the energy price $\lambda$ is not just a value but a function depending on the energy allocation $a$. Let $\lambda(\cdot)$ denote the pricing function. Then the respective payoff functions for individual agents and the coordinators become
\begin{equation}
\label{agentpayoffreversestackel}
U_i(a,\lambda(a);\theta_i), \quad \forall i=1,\ldots,N,
\end{equation}
where $a_i\in \mathcal{A}_i$ and $\lambda(a)\in\Lambda$, and
\begin{equation}
\label{coordintorpayoffreversestackel}
U_0(a,\lambda(a);\theta),
\end{equation}
where $a\in\mathcal{A}_s$ and $\lambda(a)\in\Lambda$. We assume that the coordinator knows all the private information, and announces a pricing function $\lambda(\cdot)$ to all the resource agents. Each resource agent subsequently chooses $a_i$ to maximize its payoff function after receiving $\lambda(\cdot)$. Formally, this can be captured by the following transactive energy system,
\begin{itemize}
\item agent preference: the resource agents and the coordinator agent have preferences defined by (\ref{agentpayoffreversestackel}) and (\ref{coordintorpayoffreversestackel}), respectively.
\item control decision: $\gamma_0=\lambda(\cdot)$ and $\gamma_i=a_i, \forall i=1,\ldots,N$.
\item information structure: type-dependence graph as in Figure \ref{privateinfographexample1}, and decision-dependence graph as in Figure \ref{decisiongraphexample2}.
\item solution concept: given $\lambda(\cdot)$, the decisions of resource agents are coupled, and the solution concept is the Nash equilibrium.
\end{itemize}

Given the pricing function $\lambda(\cdot)$ for energy, the problem of the second-stage is a game due to the coupling among the payoff functions of resource agents. We define the game equilibrium as follows:
\begin{equation*}
\label{lowerlvelreverse}
a_i^*(\lambda(\cdot))=\argmax_{a_i\in A_i} U_i(a_i,a_{-i}^*,\lambda(a_i,a_{-i}^*);\theta_i), \quad \forall i=1,\ldots,N,
\end{equation*}
where $a_i^*(\lambda(\cdot))$ is the optimal response of each resource agent to the pricing function $\lambda(\cdot)$ of energy. With slight abuse of notation, we denote $a^*=(a_1^*,\ldots, a_N^*)$ as $a^*(\lambda(\cdot))$ to explicitly capture the dependence of $a^*$ on $\lambda(\cdot)$. Using this notation, the decision of the coordinator agent reduces to the solution to the following optimization problem,
\begin{IEEEeqnarray}{r'l}
\label{reversestackelberggame}
\max_{\lambda(\cdot)\in\Lambda} & U_0(a^*(\lambda(\cdot)),\lambda(a^*(\lambda(\cdot)));\theta)  \IEEEyesnumber \IEEEyessubnumber \label{reversestackelberggamea} \\
\text{s.t} & a^*(\lambda(\cdot))=\argmax_{a_i\in A_i} U_i(a_i,a_{-i}^*,\lambda(a_i,a_{-i}^*);\theta_i). \IEEEyessubnumber \label{reversestackelberggameb}
\end{IEEEeqnarray}

\begin{example}
Consider again the motivating example of supplying the electricity to a group of $N$ distribution loads in the presence of constrained transmission capacity.
Different from previous examples, the payoff functions of resource agents and the coordinator agent are modeled as (\ref{agentpayoffreversestackel}) and (\ref{coordintorpayoffreversestackel}), respectively. We assume that the coordinator agent knows all the private information $\theta$. This coordinator determines the pricing function $\bar{\lambda}(\cdot)$ of energy in the first stage, and each agent determines $a_i$ accordingly in the second stage. Given the pricing function $\bar{\lambda}(\cdot)$, the optimal responses of resource agents are given as
\begin{equation*}
\label{agentresponse2}
a^*(\bar{\lambda}(\cdot))=\argmax_{a_i\in\mathcal{A}_i} V_i(a_i;\theta_i)-\bar{\lambda}(a_i,a_{-i}^*)a_i.
\end{equation*}
Taking into account the optimal responses of resource agents, the coordinator agent can determine the optimal pricing function $\bar{\lambda}(\cdot)$ by formulating the following optimization problem,
\begin{IEEEeqnarray}{r'l}
\max_{\bar{\lambda}(\cdot)\in\mathbb{R}} & U_0\left(a^*(\bar{\lambda}(\cdot)),\bar{\lambda}(a^*(\bar{\lambda}(\cdot)));\theta\right) \nonumber \\
\text{s.t} & \sum_{i=1}^N a^*_i(\bar{\lambda}(\cdot))\leq D \nonumber \\
& a^*_i(\bar{\lambda}(\cdot))=\argmax_{a_i\geq 0} V_i(a_i;\theta_i)-\bar{\lambda}(a_i,a_{-i}^*)a_i, \quad \forall i, \nonumber
\end{IEEEeqnarray}
where $a^*(\bar{\lambda}(\cdot))=(a_1^*(\bar{\lambda}(\cdot)),\ldots,a_N^*(\bar{\lambda}(\cdot)))$. This is a a special case of (\ref{reversestackelberggame}).
\end{example}

\subsection{State of the Art and Challenges}
In the literature,  problem (\ref{reversestackelberggame}) is often referred to as the reverse Stackelberg game \cite{groot2012reverse}, \cite{groot2014systematic}. It is also known as inverse Stackelberg game \cite{olsder2009phenomena}, \cite{averboukh2014inverse}, or as a Stackelberg game with an incentive information structure \cite{chang1981incentive}, \cite{ho1982control}. Since $a^*$ depends on $\gamma(\cdot)$, the problem enters the field of composed functions \cite{aczel1989functional}, which is known to be notoriously complex. For this reason, the mainstream for solving the reverse Stackelberg game is an indirect approach \cite{groot2014systematic}, \cite{ho1982control}:   instead of directly solving (\ref{reversestackelberggame}), the coordinator agent first relaxes (\ref{reversestackelberggame}) by assuming that each resource agent supports the coordinator agent to maximize its payoff function. This gives rise to the following team problem,
\begin{IEEEeqnarray}{r'l}
\label{teamproblem}
\max_{a\in\mathbb{A},\lambda\in\Lambda} & U_0\left(a,\lambda;\theta\right) \IEEEyesnumber
\end{IEEEeqnarray}
Let $(a^\tau,\lambda^\tau)$ denote the solution to the team problem (\ref{teamproblem}). Because (\ref{teamproblem}) is equivalent to  (\ref{reversestackelberggamea}) without the constraint (\ref{reversestackelberggameb}), the optimal value of (\ref{teamproblem}) provides an upper bound for the optimal value of (\ref{reversestackelberggame}). The key idea of the indirect approach is that if we can find a pricing function $\lambda(\cdot)$  that realizes the team solution, i.e., $a^*(\lambda(\cdot))=a^\tau$ and $\lambda(a^*(\cdot))=\lambda^\tau$, then $\lambda(\cdot)$ achieves this upper bound and  solves the reverse Stackelberg game (\ref{reversestackelberggame}). Since the team solution is merely an upper bound, it is not always achievable by a pricing function. When such a pricing function exists, problem (\ref{reversestackelberggame}) is called incentive-controllable \cite{ho1982control}. Furthermore, if the optimal pricing function is linear or affine, the problem is called linearly incentive-controllable.

Fortunately, a large class of reverse Stackelberg games is shown to be incentive-controllable when the system has only one resource agent and one coordinator \cite{groot2014systematic},  \cite{zheng1982existence}, \cite{bacsar1984affine}, \cite{groot2016optimal}.

\begin{theorem}
\label{1c1areversestackelberggame}
Let $N=1$. If the set $\Omega=\{(a_1,\lambda_1)\in \mathcal{A}_1 \times \mathbb{R} |U_1(a_1,\lambda_1; \theta_1)\geq U_1(a_1^\tau, \lambda_1^\tau;\theta_1) \}$ is strictly convex, and $U_1(a_1,\lambda_1;\theta_1)$ is differentiable with $\bigtriangledown_{a_1}U_1(a_1^\tau, \lambda_1^\tau;\theta_1)\neq 0$, then  (\ref{reversestackelberggame}) is linearly incentive-controllable with the following optimal pricing function,
\begin{equation}
\label{optimalpricingfunciton}
\lambda_1(a_1)=\lambda_1^\tau-Q(a_1-a_1^\tau),
\end{equation}
where the linear operator $Q: \mathcal{A}_1 \rightarrow \mathbb{R}$ is chosen according to $Q \nabla_{a_1}U_1(a_1^\tau, \lambda^\tau;\theta_1)=\nabla_{\lambda_1}U_1(a_1^\tau, \lambda_1^\tau;\theta_1)$.
\end{theorem}

The proof of Theorem \ref{1c1areversestackelberggame} can be found in \cite{zheng1982existence}. It implies that if the payoff function of the single resource agent is concave and continuously differentiable, and the payoff function of the coordinator agent is also continuously differentiable, then the market-based coordination problem is linearly incentive-controllable. In an effort to polish this result, a necessary and sufficient condition has been proposed in \cite{groot2016optimal}, which further relaxes the differentiability requirement and allows the payoff functions to be non-convex. In addition, nonlinear optimal pricing function was also considered in \cite{groot2014systematic}, which enlarges the class of incentive-controllable problems.

However, all these results are limited to the problems involving one coordinator and one resource. When there are multiple resource agents, the lower level decision making problem is typically coupled, thus the aforementioned results are hard to generalize. To circumvent this difficulty, one way is to design a set of pricing functions so that the resulting payoff functions of individual resource agents are identical  \cite{moulin2010efficient}, \cite{barreto2013design}, \cite{barreto2014incentives},  \cite{olsder2009phenomena}, \cite{jing1988solution}, \cite{tanaka2012dynamic}. This approach is possible if the coordinator has enough freedom to manipulate the system (e.g., the coordinator is allowed to design different pricing functions for different resource agents). In this case, the multi-agent problem reduces to a one-agent problem, for which numerous tools are available in the literature. However, in many transactive energy systems, tailoring the payoff functions of resource agents is not possible. In these cases, the reverse Stackelberg game is only proved to be linearly incentive-controllable when the payoff functions of resource agents are linear and the payoff function of the coordinator agent is quadratic \cite{gharesifard2014designing}. Beyond the linear-quadratic case, some preliminary attempts have also been made. For instance,  \cite{averboukh2014inverse} proves the existence of an optimal solution, and \cite{chen2012design} and  \cite{samadi2014real} derived sub-optimal solutions for the reverse Stackelberg game. However, solving the reveres Stackelberg game with multiple followers to the global optimality remains an interesting open problem \cite{groot2012reverse}.

\section{Mechanism Design}
The results on the Stackelberg game and reverse Stackelberg game build upon an underlying assumption that the coordinator agent has the complete information regarding the preferences of all the resource agents. For instance, to implement the optimal policy in Theorem \ref{1c1areversestackelberggame}, the coordinator agent has to know the gradient of the payoff function of each resource agent, which is often unavailable in practice. This section addresses this concern by formulating a few transactive energy systems that can be cast as a mechanism design problem \cite{diamantaras2009toolbox}, \cite{borgers2015introduction}. We study different results by varying the solution concepts, and survey the existing results in each category.

\subsection{Preliminaries on Mechanism Design}
Consider again the market-based coordination problem for a group of $N$ distribution loads. When it is cast as a mechanism design problem, the coordinator agent designs what information to be requested from resource agents and how to use that information so that a system-level coordination objective can be achieved.
We define an outcome $x\in \mathcal{X}$ that consists of the energy allocations $a_i$ for resource agents and their payments $t_i$ received from the coordinator agent, i.e., $x=(a_1,\ldots,a_N,t_1,\ldots,t_N)$, where $t_i=\lambda_ia_i$. Note that for resource agents that consume the energy, the corresponding payments will be negative. With slight abuse of notation, we denote the payoff function of each resource agent as $U_i(x;\theta_i)$, $\forall i=1,\ldots,N$. Let $f(\theta):\Theta \rightarrow \mathcal{X}$ denote the social choice function. The value of $f(\theta)$ specifies the outcome preferred by the coordinator when the private information is known. Therefore, the payoff function of the coordinator agent is encoded in the
social choice function. 

To elicit private information, the coordinator agent requests each resource agent to submit a bid. This is modeled by a collection of message spaces $M=M_1\times \cdots \times M_N$, where each resource agent chooses a bid $m_i\in M_i$. In addition, the coordinator also specifies  an outcome function $g(m_1,\ldots,m_N):M \rightarrow \mathcal{X}$ that determines the market outcome $x$ based on the bids from all the resource agents. This outcome function and the message space constitute a {\em mechanism}, denoted as $\Gamma=(M,g(\cdot))$. The mechanism couples the decisions of all the resource agents and induces a game problem. When different mechanisms are used, different games will be induced.


If there is a mechanism that achieves the objective of the coordinator, we say it can implement the social choice function \cite{maskin2002implementation}, \cite{jackson2001crash}. Formally, the concept of {\em implementation} is defined as follows.
\begin{definition}[\text{Implementation} \cite{mas1995microeconomic}]
\label{implementation}
A mechanism $\Gamma=(M,g(\cdot))$ implements the social choice function $f(\cdot)$ in a dominant strategy equilibrium (Bayesian Nash equilibrium, or Nash equilibrium) if for all $\theta\in \Theta$, there exists a dominant strategy equilibrium (Bayesian Nash equilibrium, or Nash equilibrium) $m^*$ of $\Gamma$ such that $g(m_1^*,\ldots,m_N^*)=f(\theta)$.
\end{definition}

In principle, one shall search over the space of all mechanisms to check whether a social choice function can be implemented.  This is clearly a daunting task. Fortunately, due to the celebrated revelation principle \cite{gibbard1973manipulation}, \cite{austen1999positive}, it turns out that it is sufficient to search over {\em direct mechanisms}.
\begin{definition}
A mechanism $\Gamma=(M,g(\cdot))$ is a direct mechanism if the message space for each resource agent is restricted to $\Theta_i$, i.e., $M_i=\Theta_i$.
\end{definition}

In a direct mechanism, the coordinator agent directly asks individual resource agents to submit their types and sets the outcome function to be the social choice function. Clearly, the coordinator expects individual resources to truthfully submitm their types. If truthful bidding is a game equilibrium, the mechanism is {\em incentive compatible} in the solution concept of dominant strategy equilibrium, Bayesian Nash equilibrium or Nash equilibrium. Equivalently, the mechanism truthfully implements the social choice function in the corresponding solution concept. According to the revelation principle, if there exists a mechanism that implements the social choice function in a solution concept, then there is also a direct mechanism that truthfully implements the social choice function in this solution concept. Therefore, we can restrict our attention to the class of direct mechanisms without loss of generality. In the rest of this section, we study how to implement a social choice function under different solution concepts, and we focus on direct mechanisms. 

\subsection{Dominant Strategy Implementation}
Consider a transactive energy system with the following formulation,
\begin{itemize}
\item agent preference: the preference of resource agent $i$ is captured by $U_i(x;\theta_i)$ with $x\in \mathcal{X}$, and the preference of the coordinator is encoded in $f(\cdot)$.
\item control decision: $\gamma_i=m_i, \forall i=1,\ldots,N$ and $\gamma_0=g(\cdot)$.
\item information structure: type-dependence graph as in Figure \ref{decisiongraphexample1}, and decision-dependence graph as in Figure \ref{decisiongraphexample2}.
\item solution concept: given $g(\cdot)$, the solution concept is the dominant strategy equilibrium.
\end{itemize}

This transactive energy system solves a mechanism design problem with the dominant strategy equilibrium, which has been the focus of many works \cite{mas1995microeconomic}, \cite{diamantaras2009toolbox}.
It turns out that the  solution is closely related to the dictatorship of the social choice function.
\begin{definition}
\label{definitionofdictorship}
A social choice function $f(\cdot)$ is dictatorial if there is an agent  $i\in \{1,\ldots,N\}$ such that for all $\theta\in \Theta$ 
\begin{equation*}
f(\theta)\in \left\{x\in \mathcal{X}: U_i(x;\theta_i)\geq U_i(x';\theta_i) \text{ for all } x'\in \mathcal{X} \right\}.
\end{equation*}
\end{definition}
\vspace{0.2cm}
Based on Definition \ref{definitionofdictorship}, if a social choice function is dictatorial, then there is an agent in the system (also referred to as the dictator) whose optimal decision is completely aligned with the social choice. In this case, the objective of the coordinator agent is equivalent to maximizing $U_i$, which totally neglects the payoff of other agents in the system. In this regard, dictatorship is an undesirable property for the social choice function. However, as dominant strategy equilibrium is rather strong, it turns out that a social choice function has to be dictatorial in order to be truthfully implementable in a dominant strategy equilibrium. This negative result is known as the Gibbard-Satterthwaite Theorem \cite{mas1995microeconomic}, \cite{barbera1990strategy}.
\begin{theorem}
\label{GStheorem}
Let $U_i(x;\theta_i)$ be a continuous functions of $x$, and assume that the range of a social choice function $f(\cdot)$ has at least three elements. Then the social choice function $f(\cdot)$ is truthfully implementable in a dominant strategy equilibrium if and only if it is dictatorial.
\end{theorem}

This theorem implies that it is impossible to find a solution to the transactive energy system described in this section, unless the social choice function is dictatorial. In practice, for the problem of market-based coordination for DERs, the objective of the coordinator agent is often to maximize the social welfare, or the profit of the coordinator. In general, the social choice function corresponding to such an objective cannot admit a dictator in the system. Hence, finding a solution of dominant strategy equilibrium to this coordination problem is generally impossible.

To circumvent this difficulty, we can either consider a more restricted class of payoff functions, or relax the solution concept. In the rest of this subsection, we consider the first case.

Consider a special case for which each resource agent has a quasi-linear payoff function denoted as $U_i(x;\theta_i)=V_i(a;\theta_i)+t_i$. The first term $V_i$ is the utility function, which is assumed to depend on $a$ instead of only $a_i$ as in (\ref{CEagents}), and the second term $t_i$ is the payment from the coordinator agent to the resource agent. Note that In addition, we assume that the objective of the coordinator is to maximize the social welfare. In other words, if we decompose the social choice function in two parts, i.e., $f(\theta)=(a^*(\theta),t^*(\theta))$, then the first part should satisfy the following, 
\begin{equation}
\label{allocativeefficiency}
a^*(\theta)=\argmax_{a\in \mathcal{A}_s } \sum_{i=1}^N V_i(a;\theta_i),
\end{equation}
Due to the special structure of this payoff functions, we can find a mechanism that truthfully implements the social choice function in a dominant strategy equilibrium. This mechanism is known as the Vickrey-Clarke-Groves (VCG) mechanism \cite{clarke1971multipart}, \cite{groves1973incentives}.
\begin{theorem}
\label{vcgmechanism}
The social choice function $f(\theta)=(a^*(\theta),t_1(\theta),\ldots,t_N(\theta))$ is truthfully implementable in a dominant strategy equilibrium if $a^*(\theta)$ satisfies (\ref{allocativeefficiency}), and
\begin{equation}
t_i(\tilde{\theta})=\sum_{j\neq i} V_j(a^*(\tilde{\theta});\tilde{\theta}_j)+h_i(\tilde{\theta}_{-i}).
\label{vcgpayment}
\end{equation}
for all $i$ and all $\tilde{\theta}_i\in \Theta_i$, where $\tilde{\theta}=(\tilde{\theta}_1,\ldots,\tilde{\theta}_N)$ is the agent bid, and $h_{-i}(\cdot)$ can be any function that does not depend on $\tilde{\theta}_i$. 
\end{theorem}
The proof of this theorem can be found in \cite[Chap 23]{mas1995microeconomic}. It implies that if the mechanism adopts the allocation rule (\ref{allocativeefficiency}) and payment rule  (\ref{vcgpayment}), then each resource agent will bid truthfully regardless of the decisions of other agents.  In fact, we can also prove that under mild conditions, the VCG mechanism  (\ref{allocativeefficiency}) and (\ref{vcgpayment}) is the only mechanism that implements the social choice function in a dominant strategy equilibrium. This uniqueness result is discussed in \cite[P.879]{mas1995microeconomic}.

The VCG mechanism has been successfully applied in the design of electricity market and demand response problems \cite{samadi2012advanced}, \cite{nekouei2015game}, \cite{samadi2011optimal},  \cite{de2016market}, \cite{xu2017efficient}, \cite{ma2018strategic}. In the design of electricity markets, the coordinator represents the independent system operator (ISO), and the resource agents are either consumers or suppliers. Since the ISO is non-profit, it is important to ensure that the mechanism can be implemented by the ISO without any external funding, i.e., $\sum_{i=1}^N t_i\leq 0$. Formally, a mechanism is {\em budget balanced} (or  weakly budget balanced) if $\sum_{i=1}^N t_i=0$ (or $\sum_{i=1}^N t_i\leq 0$). Unfortunately, the VCG mechanism fails this condition \cite{green1977characterization}. Furthermore, it is shown that under mild conditions, no social choice function can be truthfully implementable in a dominant strategy equilibrium with a budget balanced mechanism \cite{green1979incentives}. 

\subsection{Bayesian Nash Implementation}
To address the budget concern, one way is to relax the solution concept. In this subsection, we consider the transactive energy system solving a mechanism design problem with the Bayesian Nash equilibrium. We assume that the types of individual resources, $\theta=(\theta_1,\ldots,\theta_N)$, are drawn from a set $\Theta$ according to the probability distribution $\phi(\theta)$. Given the joint distribution, resource agents choose their control decisions to maximize their expected payoff. Hence, this tranactive energy system has the following formulation,
\begin{itemize}
\item agent preference: the preference of each resource agent is captured by $U_i(x;\theta_i)$ with $x\in \mathcal{X}$, and the preference of the coordinator agent is encoded in $f(\cdot)$.
\item control decision: $\gamma_i=m_i, \forall i=1,\ldots,N$ and $\gamma_0=g(\cdot)$.
\item information structure: $\phi(\theta)$ is common knowledge, type-dependence graph as in Figure \ref{decisiongraphexample1}, and decision-dependence graph as in Figure \ref{decisiongraphexample2}.
\item solution concept: given $g(\cdot)$, the decisions of resource agents are coupled, and the solution concept is the Bayesian Nash equilibrium.
\end{itemize}

Since the Bayesian Nash equilibrium is weaker than the dominant strategy equilibrium, we can slightly modify the VCG mechanism (\ref{vcgpayment}) to provide a solution that truthfully implements the social choice function in a Bayesian Nash equilibrium. Furthermore, under a proper choice of $h_i(\cdot)$ in (\ref{vcgpayment}), the mechanism can be budget balanced \cite{d1979incentives}, \cite{arrow1979property}.
\begin{theorem}
\label{vcgmechanism4bayesian}
The social choice function $f(\theta)=(a^*(\theta),t_1(\theta),\ldots,t_N(\theta))$ is truthfully implementable in a Bayesian Nash equilibrium, if $a^*(\theta)$ satisfies (\ref{allocativeefficiency}), and
\begin{align}
\label{vcgpayment4bayesian}
t_i(\tilde{\theta})= &\sum_{j\neq i} \mathbb{E}_{\tilde{\theta}_{-i}} \left[ V_j(a^*(\tilde{\theta});\tilde{\theta}_j) \left|\vphantom{\tilde{\theta}}\right.\theta_i\right] \\
&+\dfrac{1}{N-1}\sum_{k\neq i} \mathbb{E}_{\tilde{\theta}_{-i}} \left[\sum_{j\neq k} V_j(a^*(\tilde{\theta});\tilde{\theta}_j) \left|\vphantom{\sum_{k\neq i}}\right.\theta_k  \right].
\end{align}
Furthermore, the mechanism is budget balanced, i.e., $\sum_{i=1}^N t_i(\theta)\equiv 0$.
\end{theorem}
This result implies that if the coordinator determines the allocation according to (\ref{allocativeefficiency}) and set the payment based on (\ref{vcgpayment4bayesian}), then the mechanism is budget balanced, and can truthfully implement the social choice function in a Bayesian Nash equilibrium.  This holds as long as  the payoff functions of resource agents are quasi-linear, i.e., $U_i(x;\theta_i)=V_i(a;\theta_i)+t_i$. Stronger results can be obtained in a more restricted setting. For instance, if the resource agent has a scalar-type linear utility, i.e., $V_i(a;\theta_i)=\theta_iv_i(a)$, necessary and sufficient conditions are available to identify which social choice function is truthfully implementable \cite[Chap. 3]{borgers2015introduction}, \cite{cao2012optimal}.

In many application of interest, the resource agents are free to opt out of the mechanism. Therefore, in addition to guarantee incentive compatibility and budget balance, we also need to ensure {\em individual rationality constraint}, so that each resource agent wants to participate into the mechanism. Formally, the individual rationality constraint is defined as follows
\begin{definition}
\label{individualrational}
A mechanism with the outcome function $g(\cdot)$ is interim individual rational if for all $i$, we have $\mathbb{E}_{\tilde{\theta}_{-i}} [ U_i(g(\theta_i,\tilde{\theta}_{-i});\theta_i) | \theta_i ] \geq 0$ for all $\theta_i\in \Theta_i$.
\end{definition}
If a mechanism is interim individual rational, then each resource agent is willing to participate after he observes his private value and before the mechanism outcome is revealed. While individual rationality is a desirable property, according to the Myerson-Stterthwaite Theorem \cite{myerson1983efficient}, there is no mechanism that is both Bayesian Nash incentive compatible, budget balanced, and interim individual rational. This result is true even if there are only two agents holding private valuations toward a single item.

\subsection{Nash Implementation}
To achieve both individual rationality and budget balance, we further relax the solution concept to the Nash equilibrium. As such, consider the transactive energy system with the following formulation,
\begin{itemize}
\item agent preference: the preference of each resource agent is captured by $U_i(x;\theta_i)$ with $x\in \mathcal{X}$, and the preference of the coordinator agent is encoded in $f(\cdot)$.
\item control decision: $\gamma_i=m_i, \forall i=1,\ldots,N$ and $\gamma_0=g(\cdot)$.
\item information structure: type-dependence graph as in Figure \ref{decisiongraphexample1}, and decision-dependence graph as in Figure \ref{decisiongraphexample2}.
\item solution concept: adopt the Nash equilibrium as the solution.
\end{itemize}

It is possible to design  a mechanism that is budget-balanced, individual rational and  truthfully implement the desired social choice function in a Nash equilibrium. Below we discuss some examples. 

First, consider a case with the quasi-linear payoff functions defined as $U_i(x;\theta_i)=V_i(a_i;\theta_i)+t_i$. Assume that $a_i\in \mathbb{R}$ and $V_i(\cdot;\theta_i)$ is concave and differentiable. One relevant  mechanism that realizes efficient energy allocation is the scalar strategy VCG (SSVCG) mechanism \cite[Chap. 21]{nisan2007algorithmic}, \cite{johari2009efficiency}, \cite{yang2006vcg}. Similar to the VCG mechanism, the idea of SSVCG mechanism is to ask each reousrce agent to submit a utility function. However, instead of an arbitrary class of functions, the scalar strategy VCG has a more restricted message space so that the agent has to choose from a given family of functions that can be parameterized by a scalar $\sigma_i\in \mathbb{R}$, i.e., $\bar{V}_i(a_i;\sigma_i)$. In other words, the coordinator specifies a class of single-parameterized functions $\bar{V}_i(a_i;\sigma_i)$ and each agent submits a parameter $\sigma_i$.
\begin{definition}
\label{assumsptionssvcg}
Let $\mathcal{V}$ denote the family of functions that satisfy the following conditions,
\begin{enumerate}
\item For each $\sigma_i$, $\bar{V}_i(\cdot;\sigma_i)$ is strictly concave, strictly increasing, and continuously differentiable at all $\sigma_i>0$;
\item For all $\gamma>0$, and $a_i\in \mathbb{R}$, there exists a $\sigma_i$ such that $\bar{V}_i'(a_i;\sigma_i)=\gamma$.
\end{enumerate}
\end{definition}
Note that a lot of structures of $\bar{V}_i$ satisfy Definition \ref{assumsptionssvcg}, such as $\bar{V}_i(a_i;\sigma_i)=\sigma_i a_i^2$ or $\bar{V}_i(a_i;\sigma_i)=\sigma_i log(a_i)$. The coordinator can choose from any one of them. After $\mathcal{V}$ is announced, each resource agent submits a $\sigma_i$, and the coordinator runs the VCG mechanism on the payoff functions $\bar{V}_i(a_i;\sigma_i)$. It can be shown that the Nash equilibrium of the game is efficient.
\begin{theorem}
\label{ssvcgmechanism}
Consider a mechanism $\Gamma=(\mathcal{V},\ldots,\mathcal{V}, g(\cdot))$, where the outcome function $g(\cdot)$ consists of an allocation rule $a^*(\cdot)$ that satisfies
\begin{equation}
a^*(\sigma)=\argmax_a \sum_{i=1}^N \bar{V}_i(a_i;\sigma_i),
\end{equation}
and a payment rule $t_i(\cdot)$ that satisfies:
\begin{equation}
t_i(\sigma)=\sum_{j\neq i}\bar{V}_j(a^*_j(\sigma);\sigma_j)+h_i(\sigma_{-i}).
\end{equation}
Under Assumption \ref{assumsptionssvcg}, the mechanism $\Gamma$ has a Nash equilibrium $\sigma^*$ that satisfies $\bar{V}_i'(a^*(\sigma^*))=U_i(a^*(\sigma^*);\theta_i)$, and it is efficient, i.e.,  $\sigma^*=\argmax_a \sum_{i=1}^N U_i(a_i;\theta_i)$.
\end{theorem}
The proof of Theorem \ref{ssvcgmechanism} can be found in  \cite{johari2009efficiency} and \cite[Chap. 21]{nisan2007algorithmic}. It implies that at a Nash equilibrium, each resource agent chooses a parameter so that the declared marginal utility equals the true marginal utility. Note that this mechanism provides very great flexibility and has very light communication and computation requirements. For instance, if $\bar{V}_i$ is chosen to be quadratic, then each resource agent only needs to submit a scalar, and the coordinator agent only needs to solve a quadratic program to determine the allocation and payment. However, since the resource agents do not have complete information, iterations are needed to reach the Nash equilibrium.

Many other works offer efficient mechanisms when the payoff function is quasi-linear. For instance, \cite{stoenescu2006pricing} and \cite{rasouli2014electricity} consider an electricity market design where each resource agent is required to propose an allocation and price. The coordinator then determines the allocation and price so that the induced Nash equilibrium maximizes the social welfare. In \cite{maheswaran2006efficient} and \cite{maheswaran2004social},  an infinite subclass of efficient mechanisms was obtained to always maximize the social welfare, and one can choose among this class to optimize other metrics of the mechanism.
In addition, the objective of some applications is not to achieve efficient allocation,  but rather to elicit private information for other coordination purposes. For instance, \cite{muthirayan2016mechanism} considers the problem of eliciting a baseline to meet random load reduction, and \cite{mhanna2014towards} and \cite{mhanna2016faithful} elicit parameters for the demand response assets to lower payments and better predict day-ahead energy consumption. Note that all these mechanisms can be made budget balanced \cite{rasouli2014electricity}, \cite{mhanna2014towards}, \cite{mhanna2016faithful} and/or individual rational \cite{rasouli2014electricity}, \cite{stoenescu2006pricing},  \cite{muthirayan2016mechanism}.

Another notable class of mechanisms considers a uniform price market relying on the supply function bidding \cite{klemperer1989supply}. In supply function bidding, each supplier declares the amount they would like to provide at any price, and the market is cleared at a uniform price. When each resource agent becomes strategic, the supply function bidding leads to inefficient energy allocation, but has bounded efficiency loss \cite{johari2011parameterized}, \cite{johari2004efficiency}. This result was later extended to electricity market design with capacity constraints \cite{xu2016demand}, or distribution network constraints \cite{xiao2016supply}, \cite{lin2016parameterized}.

\section{Conclusion}
This paper presented a unifying framework for transactive energy systems that are designed to solve the problem of market-based coordinationa for DERs. The proposed framework consists of four elements that are important in identifying and comparing different transactive energy systems in the literature. We used these elements to analyze a number of important classes of transactive energy systems, including the competitive equilibrium, Stackelberg game, reverse Stackelberg game, and mechanism design, each of which has a special formulation. We also surveyed available tools and results in the literature. Future work includes extending the framework to explicitly capture more complicated dynamics, and incorporating uncertainties from the model and the environment.

\section{Acknowledgements}
The first two authors contribute equally to this work. This work was partly supported by the National Science Foundation under Grant CNS-1552838.


\bibliographystyle{unsrt}
\bibliography{marketbasedapproach}

\end{document}